\documentclass[a4paper,twoside,10pt,reqno]{article}
\usepackage[english]{babel}
\usepackage[dvips]{graphicx}
\usepackage[T1]{fontenc}
\usepackage[latin1]{inputenc}
\usepackage{amssymb,amsthm,dsfont,mathrsfs}
\usepackage{amsfonts,amsmath,euscript}
\usepackage{color,cite}

\newtheorem{defi}{Definition}[section]
\newtheorem{teo}{Theorem}[section]
\newtheorem{pro}[teo]{Proposition}

\newtheorem{ex}{Example}[section]

\newcounter{example}[section]

           \addto{\captionsenglish}{}

\newcommand{\hs}{\hspace{3pt}}

\newcommand{\dem}{{\bf Proof. }}
\newcommand{\fdem}{$\square$}

\newcommand{\titulo}[1]{\mbox{} \\ \noindent \textit{\textbf{\Large #1}}\\}

\renewcommand{\abstract}[1]{{\small \noindent \textbf{Abstract:} #1\\}}

\newcommand{\pchave}[1]{{\small \noindent \textbf{Keywords:} #1\\}}

\begin{document}

\begin{center}
\titulo{Dependence of maxima in space}
\end{center}

\vspace{0.5cm}
\begin{center}
\textbf{ Helena Ferreira (helenaf@ubi.pt)\hspace{3cm}Lu\'{\i}sa Pereira (lpereira@ubi.pt)}\\
\end{center}

\begin{center}
 Department of Mathematics, University of Beira
Interior, Covilh\~a, Portugal
\end{center}

\vspace{0.7cm}

\abstract{We propose a coefficient that measures the dependence among large values for spatial processes of maxima. Its main properties are: a) $k$ locations can be taken into account; b) it takes values in $[0,1]$ and higher values indicate stronger dependence; c) it is independent of the univariate marginal distributions of the random field; d) it can be related with the tail dependence and the extremal coefficients; e) it agrees with the concordance property for multivariate distributions; f) it has as a particular case the variogram from geostatistics; g) it can be easily  estimated.}

\pchave{random fields, extreme value theory, dependence coefficients}

\section{Definition of the multivariate variogram}

Natural models for spatial extremes are max-stable processes. They can be, for instance, good approximations for annual maxima of daily spatial rainfall and have been widely applied to real data in environmental, atmospheric and geological sciences.\\

Here $\{Z({\bold x}), {\bold x} \in \mathbb{R}^d\}$ denotes a strongly stationary random  field of maxima. That means that, for each choice of ${\bold x}_1, ...,{\bold x}_k$, the distribution of $(Z({\bold x}_1),...,Z({\bold x}_k))$ is a multivariate extreme value distribution $G_{{\bold x}_1, ...,{\bold x}_k}$ (Resnick, 1987). Its tail dependence function 
$$l_{{\bold x}_1, ...,{\bold x}_k}(t_1,...,t_k)={\small \displaystyle\lim_{u\downarrow 0}\frac{1-G_{{\bold x}_1, ...,{\bold x}_k}\left(G_{{\bold x}_1}^{-1}(1-ut_1),...,G_{{\bold x}_k}^{-1}(1-ut_k)\right)}{u}},$$
$(t_1,...,t_k)\in \mathbb{R}_{+}^k$,
characterizes fully the dependence among its marginals distributions $G_{{\bold x}_j}$, $j=1,...,k$.\\

Several dependence coefficients have been considered in order to resume the dependence among the marginals of $G_{{\bold x}_1, ...,{\bold x}_k}$: extremal coefficients, tail dependence coefficients and madogram.

{ \noindent}Next definition introduces our proposal and then we explore its  advantages and relations with the previous coefficients.\\

\begin{defi} The variogram of $(Z({\bold x}_1),...,Z({\bold x}_k))$ is the coefficient 
$$v({\bold x}_1, ...,{\bold x}_k)=1-\frac{k+1}{k-1}E\left(  \displaystyle\bigvee _{j=1}^kG_{{\bold x}_j}(Z({\bold x}_j))-\displaystyle\bigwedge _{j=1}^kG_{{\bold x}_j}(Z({\bold x}_j))\right).$$
\end{defi}

{ \noindent} We remark that 
$$v({\bold x}_1, ...,{\bold x}_k)=1-\frac{k+1}{k-1}E\left( \displaystyle\bigvee _{\{i,j\}\subset\{1,...,k\}} \mid G_{{\bold x}_i}(Z({\bold x}_i))- G_{{\bold x}_j}(Z({\bold x}_j))\mid \right).$$

{ \noindent} Therefore, for two locations ${\bold x}_i$ and ${\bold x}_j$, it holds $v({\bold x}_i,{\bold x}_j)=1-6\nu (\mid {\bold x}_i -{\bold x}_j \mid),$ where 
$$\nu (\mid {\bold x}_i -{\bold x}_j \mid)=\frac{1}{2}E\left( \mid G_{{\bold x}_i}(Z({\bold x}_i))- G_{{\bold x}_j}(Z({\bold x}_j))\mid \right)$$
 is the first-order variogram or madogram ( Matheron (1987), Cooley et al. (2006)).

\section{Main properties}
 
 \begin{pro} The coefficient $v({\bold x}_1, ...,{\bold x}_k)$ is independent of the univariate distributions $G_{{\bold x}_j}$, $j=1,...,k$ and it holds
$$v({\bold x}_1, ...,{\bold x}_k)=1-\frac{k+1}{k-1}\left(\frac{l_{{\bold x}}(1,...,1)}{1+l_{{\bold x}}(1,...,1)}-\sum_{\emptyset \neq I\subseteq\{1,...,k\}}(-1)^{\mid I\mid +1}\displaystyle\frac{l_{{\bold x}_I}(1,...,1)_I}{1+l_{{\bold x}_I}(1,...,1)_I}\right),$$
where ${\bold x}=({\bold x}_1, ...,{\bold x}_k)$ and ${\bold x}_I$ denotes the sub-vector of ${\bold x}$ with indices in $I$.
\end{pro}

The above proposition enhances that \\

{ \noindent}a) the tail dependence function  $l_{{\bold x}_1, ...,{\bold x}_k}$ of $G_{{\bold x}_1, ...,{\bold x}_k}$ determines 
$v({\bold x}_1, ...,{\bold x}_k)$;\\
{ \noindent}b) the extremal coefficients $\epsilon_{{\bold x}_I}=l_{{\bold x}_I}(1,...,1)_I$, for $I\subseteq \{1,...,k\}$ determine 
$v({\bold x}_1, ...,{\bold x}_k)$;\\
{ \noindent} c)  for $k=2$, we obtain  $\nu (\mid {\bold x}_i -{\bold x}_j \mid)=\frac{l_{{\bold x}_i,{\bold x}_j}(1,1)}{1+l_{{\bold x}_i,{\bold x}_j}(1,1)}-\frac{1}{2}$, which is the equation $(14)$ in Cooley et al. (2006).\\

\begin{pro} If $(Z({\bold x}_1),...,Z({\bold x}_k))$ is more concordant than $(Z({\bold y}_1),...,Z({\bold y}_k))$ then $v({\bold x}_1, ...,{\bold x}_k)\geq v({\bold y}_1, ...,{\bold y}_k)$.
\end{pro} 

\begin{pro} a) If $(Z({\bold x}_1),...,Z({\bold x}_k))$ has independent margins then $v({\bold x}_1, ...,{\bold x}_k)=0$.\\
{ \noindent}b) If $(Z({\bold x}_1),...,Z({\bold x}_k))$ has totally dependent margins then $v({\bold x}_1, ...,{\bold x}_k)=1$.\\
{ \noindent}c) $0\leq v({\bold x}_1, ...,{\bold x}_k) \leq 1$.
\end{pro}

The proofs of the above propositions can be found in Ferreira (2013).

\section{Estimation}

Let $(Z^{(i)}({\bold x}_1),...,Z^{(i)}({\bold x}_k))$, $i=1,...,n$, independent copies of $(Z({\bold x}_1),...,Z({\bold x}_k))$ and ${\hat{G}}_{{\bold x}_j}$  the empirical distribution function provided by   $Z^{(i)}({\bold x}_j)$, $i=1,...,n$, $j=1,...k$. The natural estimator for the variogram $v({\bold x}_1, ...,{\bold x}_k)$ is

$${\hat{v}}({\bold x}_1, ...,{\bold x}_k)=1- \frac{k+1}{k-1}\times \frac{1}{n}\displaystyle\sum_{i=1}^n \left( \displaystyle\bigvee _{j=1}^k{\hat{G}}_{{\bold x}_j}(Z^{(i)}({\bold x}_j))-\displaystyle\bigwedge _{j=1}^k{\hat{G}}_{{\bold x}_j}(Z^{(i)}({\bold x}_j))\right).$$

The estimator is strongly consistent and its asymptotic normality can be deduced from the Theorem $6$ in Fermanian et al. (2004).\\

\section{Application}

We compute  the estimates for the variogram of the amount of tritium  (pCi/L) in drinking water, for Muscle Shoals, Scottsboro  and Montgomery, three cities in Alabama State (USA).\\

 The first two are relatively close to a nuclear  power plant in northern Alabama (Browns Ferry reactors at Decatur): Muscle Shoals is about $20$ miles west of the plant, Scottsboro is about $60$ miles east of the plant and  both are situated on the banks of the Tennessee River. Measurements (quarterly) can be accessed at 
http://oaspub.epa.gov/enviro/erams.

The tritium data in drinking water show both Muscle Shoals and Scottsboro have greater levels than Montgomery in the period 1997-2013, and we evaluate the dependence of the annual maxima at these locations throughout the estimated  variogram ${\hat{v}}$.\\

\begin{tabular}{|c|c|c|}
\hline
 location  & Scottsboro (S) & Montgomery (M)\\  \hline
 Muscle & &\\
 Shoals (MS) & ${\hat{v}}$(MS,S)=0,4 & ${\hat{v}}$(MS,M)=0,13\\  \hline
 Scottsboro &  & ${\hat{v}}$(S,M)=0,22\\  \hline
 \end{tabular}
 
  \begin{tabular}{|c|}
\hline
${\hat{v}}$(MS,S,M)=0,26\\  \hline
 \end{tabular}

\end{document}